\documentclass[11pt,reqno]{amsart}

\usepackage[utf8]{inputenc}
\usepackage{a4wide, amsfonts}
\usepackage[english]{babel}
\usepackage{xcolor}
\usepackage[T1]{fontenc}
\usepackage{amsmath, amssymb, epic,graphicx,mathrsfs,enumerate}    
\usepackage{cite}
\usepackage[numbers,sort&compress]{natbib}
\usepackage[all]{xy}
\usepackage{color}
\usepackage{comment}
\usepackage{enumitem}
\usepackage{indentfirst}
\usepackage{esint}
\usepackage{hyperref}
\usepackage{amsthm}
\usepackage{latexsym}
\usepackage{epsfig}
\usepackage{soul}
\usepackage{geometry}
\usepackage{doi}

\hypersetup{colorlinks=true, linkcolor=blue, filecolor=black, urlcolor=blue,
citecolor=blue}
\newtheorem{theorem}{Theorem}[section]

\newtheorem{lemma}[theorem]{Lemma}
\newtheorem{prop}[theorem]{Proposition}
\newtheorem{definition}{Definition}[section]
\theoremstyle{remark}
\newtheorem{obs}{Remark}
\newtheorem{norm}{Normalization}

\newtheorem{ackn}{Acknowledgement\ \!\!\!\!\!}

\newcommand{\R}{\mathbb{R}}

\newcommand{\e}{\varepsilon}

\newcommand{\la}{\lambda}
\newcommand{\s}{\sigma}
\newcommand{\HH}{{\rm H}}

\newcommand{\Om}{\Omega}
\newcommand{\Omen}{\Omega\setminus\mathrm{MAX}(u)}
\newcommand{\p}{\Psi}

\title[Serrin problem for ring-shaped domains: the $n$-dimensional case]{On the Serrin problem for ring-shaped domains:\\ the $n$-dimensional case}

\author[V. Agostiniani]{Virginia Agostiniani}
\author[C. Bernardini]{Chiara Bernardini}
\author[S. Borghini]{Stefano Borghini}
\author[L. Mazzieri]{Lorenzo Mazzieri}

\address[V.\,Agostiniani]{
Department of Mathematics,
University of Trento 
\newline\indent Via Sommarive, 14 - 38123 Povo (Trento)}
\email{\href{mailto:virginia.agostiniani@unitn.it}{virginia.agostiniani@unitn.it}}

\address[C.\,Bernardini]{
Department of Mathematics,
University of Rome Tor Vergata 
\newline\indent Via della Ricerca Scientifica 1, 00133 Roma}
\email{\href{mailto:bernardini@axp.mat.uniroma2.it}{bernardini@axp.mat.uniroma2.it}}

\address[S.\,Borghini]{
Department of Mathematics,
University of Naples Federico II 
\newline\indent Via Cintia, Monte S. Angelo I-80126 Napoli, Italy}
\email{\href{mailto:stefano.borghini@unina.it}{stefano.borghini@unina.it}}

\address[L.\,Mazzieri]{
Department of Mathematics,
University of Trento 
\newline\indent Via Sommarive, 14 - 38123 Povo (Trento)}
\email{\href{mailto:lorenzo.mazzieri@unitn.it}{lorenzo.mazzieri@unitn.it}}

\begin{document}

\begin{abstract}
We provide several characterizations of ring-shaped rotationally symmetric solutions to the Serrin problem in arbitrary dimensions.
\end{abstract}

\maketitle

\date{}

\tableofcontents

%----------------------------------------------------------------------

\section{Introduction}

In this work, we continue the study initiated in \cite{AgoBorMaz} on the Serrin problem for ring-shaped domains, now extending the analysis to domains in the $n$-dimensional Euclidean space. More in detail, assuming $n\ge3$, we consider pairs $(\Omega,u)$, where $\Omega\subset\R^n$ is a \emph{ring-shaped domain} with smooth boundary and $u$ is a solution to
\begin{equation}\label{1}
\begin{cases}
\Delta u=-n\qquad&\text{in }\Omega,\\
u=0&\text{on }\partial\Omega.
\end{cases}
\end{equation}
By ring-shaped domain we mean a bounded domain $\Om$ whose boundary is made up of two connected components. We will denote by $\Gamma_1$ and $\Gamma_2$ the inner and the outer connected components of $\partial\Om$, respectively, so that $\partial\Om=\Gamma_1\sqcup\Gamma_2$.

In a celebrated result, Serrin \cite{Serrin} proved that if $\Omega$ is a bounded domain in $\R^n$ which admits a solution to \eqref{1} with constant normal derivative on $\partial\Omega$, then $\Omega$ is a ball and $u$ is radially symmetric. More precisely, the solution is equivalent (in the sense of Definition \ref{def sol equiv}) to
\begin{equation}\label{serrin}
u_0(x)=\frac{1-|x|^2}{2},\,\,\qquad\text{in}\quad\Omega_0=B_1(0).  
\end{equation}
In what follows we will refer to $(\Omega_0,u_0)$ as the \textit{Serrin solution}. 
Serrin's work \cite{Serrin} has also been of outstanding importance for developing the very fruitful mathematical idea 
of moving planes, which goes back to Alexandrov \cite{Alex}, was first used in PDE theory by Serrin, and originated many fundamental results (see e.g. \cite{GiNiNi, ReichelARMA, Dam, SerrinZou, Fall} for what concerns overdetermined problems).

As Serrin's result characterizes rotationally symmetric (and radially monotone) solutions defined on a domain with a connected boundary, it is equally natural to ask whether similar characterizations hold for rotationally symmetric solutions defined on domains with a disconnected boundary (possibly not radially monotone). This is precisely the aim of the present paper. 
More precisely, we are going to take into account overdetermining conditions that force the function $u$ to be rotationally symmetric and $\Om$ to be an annulus.

In regard to this issue, using the method of moving planes, Reichel \cite{Reichel} studied the following closely related problem
\begin{equation}\label{sist reichel}
\begin{cases}
\Delta u=-n\qquad\quad&\text{in }\Omega,\\
u=a,\quad \frac{\partial u}{\partial\nu}=c_1&\text{on }\Gamma_1,\\
u=0,\quad \frac{\partial u}{\partial\nu}=c_2&\text{on }\Gamma_2.
\end{cases}
\end{equation}

More in detail, assuming $a>0$ and $0<u<a$ in $\Om$, 
%so as a consequence the constant $c_1$ is non-negative, 
Reichel \cite{Reichel} proved that $u$ is radially symmetric and $\Om$ is a circular annulus (actually, his result holds for a more general class of second-order elliptic equations). Later Sirakov \cite{Sirakov} refined Reichel's result, replacing the assumption $0<u<a$ in $\Om$ with the weaker condition $c_1\geq0$. This condition is indeed crucial to start the moving planes method. Actually, when $c_1=c_2=c$ where $c$ is a negative constant, Kamburov and Sciaraffia \cite{KamSci} proved the existence of smooth annular domains (not circular annuli) for which \eqref{sist reichel} admits a non-radial solution. Finally, we point out that a crucial feature of the method of moving planes is that proving radial symmetry provides at once monotonicity of solutions in the radial direction, so a different approach is needed if we are looking for non-monotone solutions.

In this work, the main tool is the comparison
geometry technique exploited in \cite{AgoBorMaz} to characterize rotationally symmetric solutions to problem
\eqref{1} in planar, annular domains.
This same circle of ideas was simultaneoulsy developed and employed in \cite{BoChMa} and \cite{BorMazCQG} to deal with static solutions to the Einstein equations with positive cosmological constant. Further application to the study of the critical catenoid conjecture has been proposed in \cite{Es_Ma}.
%At the same time, some of the %results we are going to 
%present can be proved by 
%suitably ``gluing'' the above %mentioned Sirakov's result.

Before proceeding, we outline some fundamental properties of the solution $(\Om,u)$ to problem \eqref{1}, while also introducing some useful notation.

Since $\Omega$ is assumed to be a bounded domain, then, by the maximum principle, $u>0$ in $\Omega$ and the function $u$ achieves its maximum in the interior of the domain. In what follows, we will use the notation
$$u_{\mathrm{max}}:=\max_{\Omega}{u}\qquad\text{and}\qquad \mathrm{MAX}(u):=\{x\in\Omega\,|\,u(x)=u_{\mathrm{max}}\}.
$$
We also recall that $u$ is a real analytic function (refer e.g. to \cite{Morrey})
and that, in turn, the following facts hold.
\begin{itemize}
\item [(F1)] First, the level sets of $u$ have locally finite $\mathcal{H}^{n-1}$-measure (refer to \cite[Theorem 6.3.3]{KraPar}) and the critical level sets are discrete (see \cite[Theorem 1]{SouSou}). In particular, since the level sets of $u$ are compact for this problem, their hypersurface area is finite.
\item [(F2)] From the {\L}ojasiewicz Structure Theorem (refer e.g. to \cite[Theorem 6.3.3]{KraPar}), we have that the set of the critical points $\mathrm{Crit}(u)$
is locally a (possibly disconnected) stratified analytic sub-variety whose strata have dimensions between $0$ and $n-1$. In particular, the set of the maxima of $u$
can be decomposed as follows
$$\mathrm{MAX}(u)=\Sigma^0\sqcup\dots\sqcup \Sigma^{n-1}$$
where $\Sigma^i$ is a finite union 
of $i$-dimensional analytic submanifolds,
for every $i=0,\dots,n-1$. The set $\Sigma^{n-1}$ 
is indeed a single closed real analytic hypersurface and is usually called the \textit{top stratum} of $\mathrm{MAX}(u)$ (see \cite[Theorem 3.3]{BoChMa}): in particular, a unit normal vector field is defined everywhere on it.
\end{itemize}

\noindent Notice that our problem is invariant by translations and rescalings, i.e. if $(\Omega,u)$ is a solution to \eqref{1} then, for every $\la>0$ and $c\in\R$ also $(\Om_{\la,c},u_{\la,c})$
is a solution to \eqref{1}, where 
\begin{equation}\label{scaling inv}
\Om_{\la,c}:=\{\la x+c\,|\,x\in\Om\}\qquad\text{and}\qquad u_{\la,c}(y):=\la^2 u\left(\frac{y-c}{\la}\right).
\end{equation}
For future convenience, we provide the following definition.

\begin{definition}[Equivalent solutions]\label{def sol equiv}
We say that $(\Om,u)$ and $(\bar{\Om},\bar{u})$ solving \eqref{1} are equivalent if there exist $\la>0$ and $c\in\R$ such that $(\bar{\Om},\bar{u})=(\Om_{\la,c},u_{\la,c})$ are related as in \eqref{scaling inv}.   
\end{definition}

Before outlining the methods and tools employed to prove our results, we set up the formalism used throughout the paper.
The value of $|\nabla u|$ at the boundary of $\Om$, which is related to the so-called \emph{Wall Shear Stress} in the fluid-dynamics literature, will have an important role in the following. In this regard, it is important to note that the value of $|\nabla u|$ on the boundary components heavily depends on 
$u_{\mathrm{max}}$. 
At the same time, it would be desirable for two equivalent solutions, although having different maximum values, to predict the same boundary value for
$|\nabla u|$.
In order words, we need a scaling invariant quantity: to get it, we normalize $|\nabla u|$ by $1/\sqrt{2u_\mathrm{max}}$.  
We thus introduce the notion of \textit{Normalized Wall Shear Stress} as follows.

\begin{definition}[Normalized Wall Shear Stress]
Let $(\Omega,u)$ be a solution to problem \eqref{1} and let $S$ be a connected component of $\partial\Omega$. We define the \textnormal{Normalized Wall Shear Stress} (NWSS in short) of $u$ on $S$ as 
$$\tau(u,S):=\frac{\max\limits_S |\nabla u|}{\sqrt{2u_{\mathrm{max}}}}.$$
It is natural to define a corresponding notion associated with a connected component N of 
$\Omega\setminus\mathrm{MAX}(u)$, 
that is
$$\tau(u,N):=\max_{S\in \mathcal{S}} \tau(u,S),$$
where $\mathcal{S}$ is the set of all the connected components of $\partial\Omega\cap \overline{N}$. 
\end{definition}

\begin{obs}
Notice that if $(\Om,u)$ is a solution to \eqref{1} and $N$ is a connected component of $\Omen$, then  $\Gamma_N:=\partial\Omega\cap\overline N$, cannot be empty. This fact, which can be proved arguing by contradiction and invoking a maximum principle (see \cite[Lemma 5.1]{BorMazCQG} for details), 
implies in particular that the definition of $\tau(u,N)$ is well-posed.
\end{obs}

Our argument is based on a \textit{comparison strategy}: given a solution $(\Omega,u)$ to problem \eqref{1}, we compare some of its analytic and geometric properties with the corresponding ones for ring-shaped model solutions, that we can compute explicitly. In particular, we want to compare the gradient of a solution $(\Om,u)$, with the gradient of a suitably chosen ring-shaped model solution $(\Om_R,u_R)$. So, first of all, let us introduce what we mean by ring-shaped model solutions.\\

\textbf{Ring-shaped model solutions.} The family of radial solutions to problem \eqref{1} over annular-type domains, namely domains of the form $\{x\in\R^n\,|\, a<|x|<b\}$ where $a,b\in\R$ and $a<b$, can be parameterized by a real value $R\in(0, R_\mathrm{max})$ where $R_\mathrm{max}:=\sqrt{\frac{n-2}{n}}$, 
that we call \textit{core radius}. In particular, up to translations and rescaling, every radial solution has the following form
\begin{equation}\label{ring sol}
u_R(|x|)=\frac{1-|x|^2}{2}-\frac{R^n}{(n-2)|x|^{n-2}}
\end{equation}
over the domain
\begin{equation}\label{ring sol 2}
\Omega_R=\{x\in\R^n\,|\,r_1<|x|< r_2\}
\end{equation}
where $r_1=r_1(R)$ and $r_2=r_2(R)$ are the two positive zeroes of the function
\begin{equation}\label{f_R}
  f_R(r):=1-r^{2}-\frac{2R^n}{(n-2)r^{n-2}}.  
\end{equation}
It follows immediately that when $R\to0^+$ we recover the above-mentioned Serrin solution \eqref{serrin}.
The functions $u_R$ are not monotonically decreasing in $|x|$, indeed they are positive in the interior of $\Omega_R$, they vanish on $\partial\Omega_R$ and achieve their maximum 
\begin{equation}\label{def umax}
(u_R)_\mathrm{max}=\frac{1}{2}\left(1-\frac{n}{n-2}R^2\right)
\end{equation}
on the set $\mathrm{MAX}(u_R)=\{x\in\R^n\,|\,\,|x|=R\}$. This set cuts the domain $\Omega_R$ in two connected components, that is 
$$\Omega_R\setminus\mathrm{MAX}(u_R)=A_{r_1,R}\sqcup A_{R,r_2}$$
where  
$$A_{r_1,R}:=\{r_1<|x|<R\}\qquad\text{and}\qquad A_{R,r_2}:=\{R<|x|<r_2\}.$$
Moreover, we will denote by 
$$\Gamma_{r_1}:=\{x\in\R^n\,|\,|x|=r_1\}\qquad \text{and} \qquad \Gamma_{r_2}:=\{x\in\R^n\,|\,|x|=r_2\}.$$
\vspace{1mm}

To implement our comparison argument, we have to associate a solution with a suitable ring-shaped model solution to compare with. Since ring-shaped solutions are parameterized by the parameter $R\in(0,R_\mathrm{max})$, this means that we have to associate to $(\Om,u)$ solution to \eqref{1} a suitable value of $R$, named \textit{expected core radius}. 
The rough idea consists of studying the slope of the given solution at the boundary, so we take advantage of the above notion of Normalized Wall Shear Stress. First, we define the \emph{Inner} and \emph{Outer NWSS}, and then we provide the notion of expected core radius.
For ring-shaped model solutions,
the Normalized Wall Shear Stress of $u_R$ on $\Gamma_{r_1}$ and $\Gamma_{r_2}$ can be computed explicitly:
\[
\tau(u_R,\Gamma_{r_i})=\frac{|u'_R|}{\sqrt{2(u_R)_\mathrm{max}}} \Bigg| _{|x|=r_i},
\qquad\text{with }\quad 
|u'_R(|x|)|=|x|\,\bigg|\frac{R^n}{|x|^n}-1\bigg|,
\qquad\text{for }i=1,2.
\]
As a consequence, we can define the following two functions.

\begin{definition}[Inner and Outer NWSS]
\phantom 2
\begin{itemize}
\item The \textit{Inner NWSS} 
$$\tau_1:(0,R_\mathrm{max}]\longrightarrow[\sqrt{n},+\infty)$$
is defined as 
$$\tau_1(R):=\tau(u_R,\Gamma_{r_1})=\begin{cases}
\frac{r_1\left(\frac{R^n}{r_1^n}-1\right)}{\sqrt{1-\frac{nR^2}{n-2}}}\qquad\text{for }\,0<R<R_\mathrm{max},\\
\sqrt{n}\hfill\text{for }\,R=R_\mathrm{max}.
\end{cases}$$
\item The \textit{Outer NWSS}
$$\tau_2:[0,R_\mathrm{max})\longrightarrow[1,\sqrt{n})$$
is defined as
$$\tau_2(R):=\tau(u_R,\Gamma_{r_2})=\begin{cases}
1\hfill\text{for }\,R=0,\\
\frac{r_2\left(1-\frac{R^n}{r_2^n}\right)}{\sqrt{1-\frac{nR^2}{n-2}}}\qquad\text{for }\,0<R<R_\mathrm{max}.
\end{cases}$$  
\end{itemize}
\end{definition}

\noindent Notice that $\tau_1$ and $\tau_2$ are both continuous, that $\tau_1$ is strictly decreasing, and that 
$\tau_2$ is strictly increasing.
Moreover, $\lim_{R\to0^+}\tau_1(R)=+\infty$ and $\lim_{R\to R_\mathrm{max}}\tau_2(R)= \sqrt{n}$.\\

\noindent We are now in the position to introduce the very key tool to implement our comparison argument: the notion of \textit{expected core radius}.

\begin{definition}[Expected core radius]
Let $(\Omega,u)$ be a solution to \eqref{1} and $S$ be a connected component of $\partial\Omega$. The expected core radius of $u$ on $S$ is defined as
$$R(u,S):=\tau_2^{-1}\big(\tau(u,S)\big)\qquad\text{if }\,\tau(u,S)\in [1,\sqrt{n})$$
and 
$$R(u,S):=\tau_1^{-1}\big(\tau(u,S)\big)\qquad\text{if }\, \tau(u,S)\in[\sqrt{n},+\infty)$$
Analogously, we can define the expected core radius on a connected component $N$ of $\Omega\setminus\mathrm{MAX}(u)$ as follows
$$R(u,N):=\tau_2^{-1}\big(\tau(u,N)\big)\qquad\text{if }\, \tau(u,N)<\sqrt{n}$$
and 
$$R(u,N):=\tau_1^{-1}\big(\tau(u,N)\big)\qquad\text{if }\, \tau(u,N)\ge\sqrt{n}$$
\end{definition}

\begin{obs}
It is immediate to observe that $R(u_R,\Gamma_{r_i})=R$ for $i=1,2$, since $\tau(u_R,\Gamma_{r_i})=\tau_i(R)$ by definition, analogously $R(u_R,A_{r_1,R})=R(u_R,A_{R,r_2})=R$. 
\end{obs}

Hence, given a solution $(\Om,u)$ to problem \eqref{1} and a connected component $N$ of $\Omen$, we can associate to $(\Om,u)$ a ring-shaped model solution $(\Om_R,u_R)$ defined as in \eqref{ring sol}-\eqref{ring sol 2} such that $R=R(u,N)$, the expected core radius of $u$. Our main goal is to \textit{compare} the gradient of $u$ with the gradient of the corresponding ring-shaped model solution $(\Om_R,u_R)$ and, in particular, to obtain a sharp and rigid upper bound for the former in terms of the latter (refer to Theorem \ref{stime gradiente} below). For a meaningful comparison, it is necessary to normalize the solution $(\Om,u)$ in the following way.

\begin{norm}
Let $(\Omega,u)$ be a solution to problem \eqref{1} and $N$ be a connected component of $\Omega\setminus\mathrm{MAX}(u)$. Let $R=R(u,N)$  be the expected core radius of $u$ on $N$, we rescale the domain and the function as in \eqref{scaling inv} in such a way that
\begin{equation}\label{normaliz}
u_\mathrm{max}=(u_R)_\mathrm{max}
\end{equation}
where $(u_R)_{\rm max}$ is defined as in \eqref{def umax}. This means that we consider an equivalent pair $u_{\la,c}$ (which we still denote by $u$)  where $\la=\sqrt{\frac{(u_R)_\mathrm{max}}{u_\mathrm{max}}}$ and $c=0$.
\end{norm}

%--------------------------------------------------------------------------

\section{Statement of the main results}

To effectively use the concept of expected core radius, we have first to prove that it is well-posed and non-negative. This is the content of the next theorem, which also provides a characterization of the \textit{Serrin solution} in terms of the expected core radius.

\begin{theorem}\label{theo c new}
Let $(\Omega,u)$ be a solution to problem \eqref{1} where $\Omega$ bounded domain with smooth boundary, and let $N$ be a connected component of $\Omega\setminus\mathrm{MAX}(u)$. Then, 
$$R(u,N)\,\,\, \text{is well-defined and non-negative}.$$ 
Moreover, $R(u,N)=0$ if and only if $(\Omega,u)$ is equivalent to the Serrin solution \eqref{serrin}.
\end{theorem}

Note that $R(u,N)$ is well-defined if and only if $\tau(u,N)\ge1$. This is indeed what is shown in the proof of Theorem \ref{theo c new}, in Section \ref{sez2}.

The following theorem, which is our main comparison result, can be seen as a natural counterpart of the previous theorem when the rigid model solution is ring-shaped. It says, in particular, that if the expected core radius of the outer boundary coincides with the one of the inner boundary, then the solution is rotationally symmetric. 

\begin{theorem}\label{theo d new}
Let $(\Omega,u)$ be a solution to problem \eqref{1} such that $\Omega$ is a ring-shaped domain and $\mathcal{H}^{n-1}(\mathrm{MAX}(u))>0$.
In particular the expected radii of $u$ on $\Gamma_1$ and on $\Gamma_2$, that is $R_1=R(u,\Gamma_1)$ and $R_2=R(u,\Gamma_2)$, are both well-defined. If $\tau(u,\Gamma_2)<\sqrt{n}$, then it holds
$$R_2\ge R_1 >0.$$
Moreover, $R_1=R_2=:R$ if and only if $(\Om,u)$ is equivalent to the ring-shaped model solution $(\Om_R,u_R)$ defined as in \eqref{ring sol}-\eqref{ring sol 2}.
\end{theorem}

The expected radii $R_1$ and $R_2$ turn out to be well-defined because in the proof it is indeed showed that $\tau(u,\Gamma_1)$ and $\tau(u,\Gamma_2)$ are both greater than or equal to $1$. 
More precisely, under the hypothesis of the previous theorem, one would expect  that
\[
\tau(u,\Gamma_1)\geq\sqrt n,
\qquad\qquad
\tau(u,\Gamma_2)<\sqrt n.
\]
Whereas the first inequality can indeed be proved and constitutes one 
of the steps of the proof of Theorem \ref{theo d new}, we have not been able to prove the second inequality so far. 
Let us also remark that the fact that $\tau(u,\Gamma_1)\geq\sqrt n$ can be easily proved  by using the Pohozaev identity in the connected component of $\Om\setminus{\rm MAX}(u)$ touched by $\Gamma_1$, under the simplifying assumption that
$\Gamma_1$ is star-shaped.

Under the hypothesis that the set of maximum points has positive $\mathcal{H}^{n-1}$-measure, we are also able to characterize ring-shaped model solutions to 
problem \eqref{1} coupled with some partially overdetermining conditions. More in detail, we have the following result.

\begin{theorem}\label{theo b new}
Let $(\Omega,u)$ be a solution to problem \eqref{1} such that $\Omega$ is a ring-shaped domain. Assume that $|\nabla u|$ is constant on either the inner or the outer boundary component of $\Om$. Then we have only two possibilities:
\begin{itemize}
    \item[i)] $(\Omega,u)$ is equivalent to a ring-shaped model solution \eqref{ring sol}-\eqref{ring sol 2}
    \item[ii)] $(\Omega,u)$ is such that $\mathcal{H}^{n-1} (\mathrm{MAX}(u))=0$. 
\end{itemize}
\end{theorem}

Our analysis also allows us to prove some a priori bounds on the area of the boundary components of a given region $N\subseteq \Omen$. Such bounds relate the area of a connected component of $$\Gamma_N:=\overline N\cap\partial\Om,$$
to the integral average of the mean curvature of $\Gamma_N$ itself (see Proposition \ref{prop bnd ar 1}). In another case, we relate the area of the smooth portion of the top-stratum intersecting $\overline N$, that is, \begin{equation}\label{sigma N}
\Sigma_N:=\overline{N}\cap\Sigma^{n-1},
\end{equation}
to the integral average of the mean curvature of $\Sigma_N$ itself (see Proposition \ref{prop bnd ar 2} below).
Finally, in the following proposition, we compare the area of $\Sigma_N$ to the area of $\Gamma_N$. 
All these estimates come with the corresponding rigidity result.

\begin{prop}\label{theo 2.4}
Let $(\Om,u)$ be a solution to problem \eqref{1} and $N$ be a connected component of $\Omen$. Let $R=R(u,N)\in[0, R_\mathrm{max})$ be the expected core radius of $u$ in the region $N$, $r_1$ and $r_2$ be the two positive zeros of the function \eqref{f_R}. Assume that Normalisation 1 holds. Finally, suppose that $\Gamma_N$ is connected and $\Sigma_N \not=\varnothing$. 
\begin{itemize}
    \item[i)] If $\tau(u,N)<\sqrt{n}$, then 
    $$\frac{|\Sigma_N|}{R^{n-1}}\,\le\,\frac{|\Gamma_N|}{r_2^{n-1}}.$$
    \item[ii)] If $\tau(u,N)>\sqrt{n}$, then 
    $$\frac{|\Sigma_N|}{R^{n-1}}\le\frac{|\Gamma_N|}{r_1^{n-1}}.$$
\end{itemize}
Moreover, if equality holds, then $(\Om,u)$ is equivalent to the ring-shaped model solution $(\Om_R,u_R)$.    
\end{prop}

%----------------------------------------------------------------------

\section{Well-posedness and non-negativity of the expected core radius}\label{sez2}

Notice that the notion of expected core radius associated with a given region $N$ connected component of $\Omega\setminus{\rm MAX}(u)$ is well-posed only if $\tau(u,N)\ge1$. Indeed, if it is not the case, $\tau(u,N)$ does not lie in the domain of definition of the function $\tau_1^{-1}$ or $\tau_2^{-1}$. This is the content of Theorem \ref{teo2.1new} below, in order to prove it, we need the following preliminary results.

\begin{prop}\label{prop2.2new}
Let $\Omega$ be a bounded domain with smooth boundary and let $u$ be a solution to problem \eqref{1}. Then,
$$\max_{\partial\Omega}\frac{|\nabla u|}{\sqrt{2u_{\mathrm{max}}}}\ge1.$$
Moreover, the equality holds if and only if $(\Omega,u)$ is equivalent to the Serrin solution \eqref{serrin}.
\end{prop}

\noindent In other words, the Serrin solution has the least possible Normalized Wall Shear Stress among all the solutions with connected boundary.

\proof
For simplicity of the proof, we can rephrase the above statement in the following equivalent way. Let $(\Omega,u)$ be a solution to problem \eqref{1} such that 
\begin{equation}\label{cond bordo <=1}
\frac{|\nabla u|}{\sqrt{2u_{\mathrm{max}}}}\le 1\qquad\text{on } \partial\Omega,   
\end{equation} 
then $(\Omega,u)$ is equivalent to the Serrin solution. 

The proof of this fact is based on a maximum principle argument and follows the ideas in the proof of \cite[Proposition 2.2]{AgoBorMaz}, so we skip the details. From the Bochner formula and the fact that $\Delta u=-n$, we get 
$$\Delta\left(|\nabla u|^2+2u\right)=2\left(|\nabla^2u|^2-\frac{(\Delta u)^2}{n}\right)\ge0.$$
Since by assumption $|\nabla u|\le\sqrt{2u_{\mathrm{max}}}$ on $\partial \Omega$, from the weak maximum principle, we get that $|\nabla u|^2+2u$ is constant and in particular
$$|\nabla u|^2+2u\equiv 2 u_{\mathrm{max}}\qquad\text{in }\Omega.$$
This implies that
$$0=\Delta\left(|\nabla u|^2+2u\right)=2\left(|\nabla^2u|^2-\frac{(\Delta u)^2}{n}\right)$$
hence $|\nabla^2u|^2=\frac{(\Delta u)^2}{n}$, from which it follows that, up to rescalings, $(\Omega,u)$ is rotationally symmetric.
To conclude the proof it is sufficient to recall that the only rotationally symmetric solution that fulfills \eqref{cond bordo <=1}
is the Serrin solution.
\endproof

We now introduce an auxiliary function defined on the set of values of $u$, whose monotonicity properties will be crucial to prove Theorem \ref{teo2.1new}

\begin{prop}\label{prop2.3new}
Let $(\Omega,u)$ be a solution to problem \eqref{1} where $\Omega$ is a bounded domain with smooth boundary, and let $N$ be a connected component of $\Omega\setminus\mathrm{MAX}(u)$ such that $\tau(u,N)\le1$. 
\begin{itemize}
\item The function $U:[0,u_\mathrm{max})\to\R$ defined as 
$$U(t):=\frac{1}{[2n(u_\mathrm{max}-t)]^{n/2}}\int\limits_{\{u=t\}\cap \overline{N}} |\nabla u|d\sigma.$$
is non-increasing. 
\item If $\mathcal{H}^{n-1}(\mathrm{MAX}(u)\cap\overline{N})>0$ then 
\begin{equation}\label{limU}
    \lim_{t\to u^-_\mathrm{max}}U(t)=+\infty.
\end{equation}
\end{itemize}
\end{prop}

\proof 
First of all, we point out that the function $U$ is well-defined, because $|\nabla u|$ is bounded in $\Omega$ and the level sets of $u$ have finite $\mathcal{H}^{n-1}$-measure (see (F1) in the Introduction). 
In what follows we adapt the proof of \cite[Proposition 2.3]{AgoBorMaz} to the present $n$-dimensional case, so we recall only the main ideas here. An application of the maximum principle yields the following gradient estimate 
$$|\nabla u|^2\le2u_\mathrm{max}-2u,\quad\text{on }N,$$
which in turn gives
$$\mathrm{div}\left(\frac{\nabla u}{[2n(u_\mathrm{max}-u)]^{n/2}}\right)\le0.$$
%= \frac{-n^2}{(2n(u_\mathrm{max}-u))^{n/2+1}}\Big(2u_\mathrm{max}-2u-|\nabla u|^2\Big)
Integration by parts on the set $\{t_1\le u\le t_2\}\cap\overline{N}$ together with the definition of $U$ gives
$$U(t_2)-U(t_1)\le0$$
which is the desired monotonicity.

Using the Łojasiewicz inequality and the compactness of $\mathrm{MAX}(u)$, we obtain that for $t$ sufficiently close to $u_\mathrm{max}$
\begin{equation}\label{U>=}
U(t)\ge\frac{C}{\left(2n(u_\mathrm{max}-t)\right)^{\frac{n}{2}-\theta}}\,\Big|\{u=t\}\cap\overline{N}\Big|,  
\end{equation}
where $C>0$ and $\theta\in\left[\frac{1}{2},1\right)$. To conclude, it is thus sufficient to prove that the hypothesis $\mathcal{H}^{n-1}(\mathrm{MAX}(u)\cap\overline{N})>0$ gives that 
\begin{equation}\label{limsup>0}
\limsup\limits_{t\to u^-_\mathrm{max}}\big|\{u=t\}\cap \overline{N}\big|\ge|\Sigma_N|.
\end{equation}
This fact can be showed by following the lines of the proof of \cite[Proposition 5.4]{BorMazCQG}.
\endproof

Finally, combining Proposition \ref{prop2.2new} with Proposition \ref{prop2.3new}, we obtain the following result, which can be seen as a localized version of Proposition \ref{prop2.2new}.

\begin{theorem}\label{teo2.1new}
Let $(\Omega,u)$ be a solution to problem \eqref{1} where $\Omega$ is a bounded domain with smooth boundary, and let $N$ be a connected component of $\Omega\setminus\mathrm{MAX}(u)$. 
Then 
$$\tau(u,N)\ge1.$$
Moreover, the equality holds if and only if $(\Omega,u)$ is equivalent to the Serrin solution \eqref{serrin}.
\end{theorem}

\proof
As before, we prove that if $\tau(u,N)\le1$ then $(\Omega,u)$ is equivalent to the Serrin solution. If $\tau(u,N)\le1$, we are under the assumptions of 
Proposition \ref{prop2.3new}, so from the monotonicity of $U$, we get that
$$\lim\limits_{t\to u^-_\mathrm{max}}U(t)\le U(0)=\frac{1}{(2n u_\mathrm{max}) ^{n/2}}\int\limits _{\partial\Omega\cap\overline{N}} |\nabla u|d\sigma\le\frac{1}{(2n u_\mathrm{max})^{\frac{n-1}{2}}}\big|\partial\Omega\cap\overline{N}\big|<+\infty.$$
It follows that $\mathcal{H}^{n-1}(\mathrm{MAX}(u)\cap\overline{N})=0$, otherwise $\lim_{t\to u^-_\mathrm{max}}U(t)=+\infty$ which would lead to a contradiction.
The fact that $\mathcal{H}^{n-1}(\mathrm{MAX}(u)\cap \overline{N})=0$ implies that $\Omega\setminus\mathrm{MAX}(u)$ is connected, and hence $\partial\Omega\cap \overline{N}=\partial\Omega$ this fact coupled with Proposition \ref{prop2.2new}, allows us to conclude.
\endproof

\proof[Proof of Theorem \ref{theo c new}] 
It follows by recalling the definition of expected core radius and exploiting Theorem \ref{teo2.1new}.
\endproof

%-----------------------------------------------------------------------------

\section{Gradient estimates}\label{sez4}

The above results tell us that, given a solution $(\Om,u)$ to problem \eqref{1} and $N$ connected component of $\Omen$, we are able to adequately associate with it a unique core radius $R=R(u,N)$. In particular, we can associate with each solution $(\Om,u)$ a unique ring-shaped model solution $(\Om_R,u_R)$: the one with $R=R(u,N)$. 

The aim of this section is to compare the gradient of the solution $(\Om,u)$ with the gradient of the model solution $(\Om_R,u_R)$. In order to do this, it is necessary — and will be assumed throughout — that Normalization 1 holds. 

Let us first introduce two auxiliary functions named \textit{pseudo-radial functions}. Observe that 
$$u_{R_{|_{[r_1,R]}}}:[r_1,R]\to [0,(u_R)_\mathrm{max}]\qquad\text{and} \qquad u_{R_{|_{[R,r_2]}}}:[R,r_2]\to [0,(u_R)_\mathrm{max}]$$
are monotonically increasing and decreasing, respectively. Hence, we can define
$$\psi_-:[0,(u_R)_\mathrm{max}]\to[r_1,R]\qquad\text{and}\qquad \psi_+:[0,(u_R)_\mathrm{max}]\to[R,r_2]$$
as the inverse of $u_{R_{|_{[r_1,R]}}}$ and $u_{R_{|_{[R,r_2]}}}$, respectively.
It is useful to compute 
\begin{equation}\label{def psi punto e due punt}
\dot{\psi}_\pm=-\frac{\psi_\pm^{n-1}}{\psi_\pm^n-R^n},\qquad\qquad\ddot{\psi}_\pm=\dot{\psi}^3_\pm\left(1+(n-1)\frac{R^n}{\psi_\pm^n}\right). 
\end{equation}
The main idea of our comparison procedure is the following. If $N$ is a region with low NWSS, then our reference model will be the model solution \eqref{ring sol} restricted to $A_{R,r_2}$, while if $N$ a region having high NWSS, then our reference model will be the model solution \eqref{ring sol} restricted to $A_{r_1,R}$. More precisely, we give the following definition.

\begin{definition}[Pseudo-radial functions]\label{def Phi}
Given $(\Om,u)$ solution to problem \eqref{1} and $N$ a connected component of $\Omen$, let $R=R(u,N)\in[0, R_\mathrm{max})$ be the expected core radius of $u$ in the region $N$, and $r_1=r_1(R)$ and $r_2=r_2(R)$ be the two positive zeros of the function \eqref{f_R}. 
\begin{itemize}
    \item[i)] If $\tau(u,N)<\sqrt{n}$, we define the pseudo-radial function $\Psi_+:\overline N\to[R,r_2]$ as
    $$\Psi_+(p):=\psi_+\left(u(p)\right).$$
    \item[ii)] If $\tau(u,N)>\sqrt{n}$, we define the pseudo-radial function $\Psi_-:\overline N\to[r_1,R]$ as
    $$\Psi_-(p):=\psi_-\left(u(p)\right).$$    
\end{itemize}    
\end{definition}

\begin{obs}
Notice that the pseudo-radial functions coincide with the radial coordinate $|x|$ if and only the solution $(\Omega,u)$ coincides with the model solution \eqref{ring sol}-\eqref{ring sol 2}.
\end{obs}

\noindent Let us introduce the notation
$$W:=|\nabla u|^2\qquad\text{and}\qquad W_R:=|u'_R|^2\circ\Psi$$
where, given $p\in\overline N$, we set $\Psi(p):=\Psi_+(p)$ if $\tau(u,N)<\sqrt{n}$, while $\Psi(p):=\Psi_-(p)$ if $\tau(u,N)>\sqrt{n}$.
Note that using \eqref{def psi punto e due punt} we can rewrite $W_R$ more explicitly as
\begin{equation}\label{W_R}
W_R=\left(\frac{\Psi^n-R^n}{\Psi^{n-1}}\right)^2.
\end{equation}

\begin{obs}\label{spegazione W WR}
Let us consider $p\in\Gamma_N:=\partial\Om\cap\overline{N}$. If $\tau(u,N)<\sqrt{n}$ then
$$\frac{|\nabla u|(p)}{\sqrt{2u_{\rm max}}}\le \tau(u,N)=\tau_2(R)=\frac{|u'_R|(r_2)}{\sqrt{2(u_R)_{\rm max}}}=\frac{|u'_R|\circ \Psi}{\sqrt{2(u_R)_{\rm max}}}.$$
If $\tau(u,N)>\sqrt{n}$ a similar chain of inequalities holds. Since we are assuming Normalisation 1, we have that $W\le W_R$ on $\partial \Om$.
\end{obs}

Our fundamental gradient estimate written in terms of the functions $W$ and $W_R$ is the content of the following theorem.

\begin{theorem}\label{stime gradiente}
Let $(\Om,u)$ be a solution to problem \eqref{1} and $N$ be a connected component of $\Omen$. Assume that the expected core radius $R=R(u,N)\in[0, R_\mathrm{max})$. Then, 
\begin{equation}\label{W<=WR}
    W\le W_R,\quad\text{on }N.
\end{equation}
Moreover, if there exists $p\in N$ such that $W(p)=W_R(p)$ then $(\Om,u)$ is equivalent to $(\Om_R,u_R)$.    
\end{theorem}

\proof
The idea of the proof is finding a suitable elliptic inequality for the function 
\[
F_\beta:=\beta(\p)(W-W_R),
\]
where $\beta(\p)$ is a positive function that will be set later, and then arguing via a Maximum Principle argument.
We compute
\begin{equation}\label{delta F}
\begin{split}
\Delta F_\beta=\beta\Delta(W&-W_R)+2\frac{\beta'}{\beta}\frac{\p^{n-1}}{R^n-\p^n}\langle\nabla F_\beta,\nabla u\rangle-n\frac{\p^{n-1}}{R^n-\p^n} \frac{\beta'}{\beta}F_\beta\\
&+\left[\left(\frac{\beta'}{\beta}\right)'-\left(\frac{\beta'}{\beta}\right)^2+\left(\frac{\p^n+(n-1)R^n}{\p(R^n-\p^n)}\right)\frac{\beta'}{\beta}\right]\frac{\p^{2n-2}}{(R^n-\p^n)^2}|\nabla u|^2 F_\beta.
\end{split}
\end{equation}
First, we estimate the term $\Delta(W-W_R)$. We compute it explicitly using the Bochner formula and the fact that $\Delta u=-n$, we get
\begin{equation}\label{delta W-WR}
\Delta(W-W_R)=2|\nabla^2u|^2-2n\left(1+(n-1)\frac{R^n}{\Psi^n}\right)-\frac{2n(n-1)R^n}{\p^2(R^n-\p^n)}|\nabla u|^2.
\end{equation}
In order to obtain an elliptic inequality for $F_\beta$, we take advantage of the following estimate from below for the term $|\nabla^2u|^2$
\begin{equation}\label{stima |nabla2u|2}
\begin{split}
|\nabla^2u|^2\ge-n\frac{R^n\p^{n-2}}{(R^n-\p^n)^2}&\langle\nabla(W-W_R),\nabla u\rangle\\
&+n\left(1+\frac{2(n-1)R^{2n}}{\p^2(R^n-\p^n)^2}|\nabla u|^2-\frac{(n-1)R^{2n}\p^{2n-4}}{(R^n-\p^n)^4}|\nabla u|^4\right), 
\end{split}
\end{equation}
which in turn can be obtained by computing explicitly the following inequality
\begin{equation}\label{disug hessiana}
\Bigg|\nabla^2u+n\frac{R^n\p^{n-2}}{(R^n-\p^n)^2}du\otimes du+\left(1-\frac{R^n\p^{n-2}}{(R^n-\p^n)^2}|\nabla u|^2\right)g_{\R^n}\Bigg|^2\ge0.  
\end{equation}
Equations \eqref{delta W-WR} and \eqref{stima |nabla2u|2} yield
\begin{equation}\label{stima delta W-WR}
\begin{split}
\Delta(W-W_R)\ge-2n&\frac{R^n\p^{n-2}}{(R^n-\p^n)^2}\langle\nabla(W-W_R),\nabla u\rangle\\
&+2n(n-1)\frac{R^n\p^{n-2}}{(R^n-\p^n)^2}\left(1-\frac{R^n\p^{n-2}}{(R^n-\p^n)^2}|\nabla u|^2\right)(W-W_R).
\end{split}
\end{equation}
In turn, equations \eqref{delta F} and \eqref{stima delta W-WR} give
\begin{align*}
\Delta F_\beta\ge\,&\frac{2\p^{n-1}}{R^n-\p^n}\left[\frac{\beta'}{\beta}-\frac{nR^n}{\p(R^n-\p^n)}\right]\langle\nabla F_\beta,\nabla u\rangle-
\frac{n\p^{n-1}}{R^n-\p^n}\left[\frac{\beta'}{\beta}-\frac{2(n-1)R^n}{\p(R^n-\p^n)}\right]F_\beta\\
&+\left[\left(\frac{\beta'}{\beta}\right)'-\left(\frac{\beta'}{\beta}\right)^2+\frac{\p^n+(3n-1)R^n}{\p(R^n-\p^n)}\frac{\beta'}{\beta}-\frac{2(n^2-n)R^{2n}}{\p^2(R^n-\p^n)^2}\right]\frac{\p^{2n-2}}{(R^n-\p^n)^2}|\nabla u|^2 F_\beta.
\end{align*}
Now, choosing 
\[
\beta(\p)=\frac{\p^{2(n-1)}}{(R^n-\p^n)^{2(\frac{n-1}n)}},
\]
we have that the coefficient of the term $F_\beta$ vanishes,
whereas the coefficient of the term $|\nabla u|^2F_\beta$
becomes positive. More precisely, with the above choice we get
\begin{equation}\label{eq ellittica fbeta}
\Delta F_\beta
-
\frac{2(n-2)\p^{n-2}}{R^n-\p^n}
\langle\nabla F_\beta,\nabla u\rangle
-2(n-1)(n+2)\frac{R^n\p^{3n-4}}{(R^n-\p^n)^4}
|\nabla u|^2F_\beta
\,\ge\,0
\qquad\text{on }N.  
\end{equation}
%We require that 
%$$\frac{\beta'}{\beta}=\frac{nR^n}{\p(R^n-\p^n)},$$
%this corresponds to setting
%$$\beta(\p)=\frac{\p^n}{|R^n-\p^n|}.$$
%We finally obtain
%$$\Delta F_\beta-\frac{nR^n\p^{n-2}}{(R^n-\p^n)^2}\left[n-2+(n+2)\frac{\p^{2n-2}}{(R^n-\p^n)^2}|\nabla u|^2\right]F_\beta\ge0\qquad\text{on }N.$$
Now we observe that, $F_\beta\le0$ on $\partial N$. Indeed, on the one hand, $W\le W_R$ on $\Gamma_N=\partial N\cap\partial \Om$ (see Remark \ref{spegazione W WR}) and $\beta$ is positive. On the other hand, $F_\beta$ goes to zero when approaching $\partial N\setminus \partial\Om=\mathrm{MAX}(u)\cap\overline{N}$.
This fact can be proved by rewriting $F_\beta$ as
$$F_\beta=\p\frac{W}{\sqrt{W_R}}-\p\sqrt{W_R}$$
by taking advantage of the limit $ \lim\limits_{x\to p, \,x\in N} \frac{W_R}{u_\mathrm{max}-u}=2n$ (see Lemma \ref{stime appendix} I) below) and of the reverse Łojasiewicz inequality.

We conclude the proof of \eqref{W<=WR} with the following argument. Suppose by contradiction that $\max_{\overline{N}} F_\beta>0$; since $\max_{\partial{N}} F_\beta\le 0$ then there exists $p\in N$ such that $F_\beta(p)>0$. In particular, $\nabla F_\beta(p)=0$ and $\Delta F_\beta(p)\le0$. 
This fact, coupled with \eqref{eq ellittica fbeta}, gives a contradiction.  We have thus obtained that 
$\max_{\overline{N}} F_\beta\leq0$, which gives \eqref{W<=WR}.

For what concerns the proof of the rigidity part, one can argue in two steps. First, from the Strong Maximum principle, one gets that $W=W_R$ in the whole $N$. Secondly, from the fact that
$|\nabla u|=\sqrt{W_R}$ where $W_R$ can be seen as a function of $u$ and from the fact that equality holds in \eqref{disug hessiana}, one obtains that the level sets of $u$ are totally umbilic, proving rotational symmetry of $u$ in $N$. 
The symmetry on the whole $\Omega$ is then straightforward in the case when $\Om\setminus{\rm MAX}(u)$ is connected. 
When instead $\Om\setminus{\rm MAX}(u)$ is disconnected, one first deduce rotational
symmetry in a connected component and 
secondly extends it to the whole of $\Om$,
by analyticity of $u$.
\endproof

In the following lemma, we fix a connected component $N$ of $\Omen$ and we collect some useful estimates for $W$ and $W_R$ in $N$ near $\mathrm{MAX}(u)$. We recall that $\Sigma_N$ is defined as in \eqref{sigma N} and we adopt the convention that the mean curvature $\HH$ of a hypersurface is defined as the trace of the second fundamental form.

\begin{lemma}\label{stime appendix}
Let $(\Om,u)$ be a solution to problem \eqref{1} and $N$ be a connected component of $\Omen$.  Assume that the expected core radius $R=R(u,N)\in[0, R_\mathrm{max})$.
\begin{itemize}
\item[I)] For every $p\in\mathrm{MAX}(u)$, 
\begin{equation}\label{lim WR/u-umax}
\lim\limits_{x\to p, \,x\in N} \frac{W_R}{u_\mathrm{max}-u}=2n.  
\end{equation}
Moreover, in a neighborhood of $p$ in $\overline N$,  the following expansion holds
\begin{equation}\label{stima W_R}
W_R=2n(u_\mathrm{max}-u)\mp\frac{4(n-1)\sqrt{2n}}{3R}(u_\mathrm{max}-u)^\frac{3}{2}+\mathcal{O}((u_\mathrm{max}-u)^2),
\end{equation}
taking the $-$ sign when $\tau(u,N)<\sqrt{n}$ and the $+$ sign when $\tau(u,N)>\sqrt{n}$.
\item[II)] Let us assume that $\Sigma_N\not=\emptyset$ and let us denote by $\HH(p)$ the mean curvature of $\Sigma_N$ at $p$ computed with respect to the unit normal pointing outside $N$, and by $r$ the signed distance from $\Sigma_N$ that is positive in $N$. Then for every $p\in\Sigma_N$, there exists $V$ neighborhood of $p$ such that $V\cap\mathrm{MAX}(u)\subset\Sigma_N$ where
\begin{align*}
&W=n^2r^2+\HH\, n^2r^3+\mathcal{O}(r^4),\\
&W_R=n^2r^2+n^2\left(\frac{\HH}{3}\mp\frac{2n-2}{3R}\right)r^3+\mathcal{O}(r^4). 
\end{align*}
Here, we take the $-$ sign when $\tau(u,N)<\sqrt{n}$ and the $+$ sign when $\tau(u,N)>\sqrt{n}$. 
\item[III)] $$\lim_{t\to u^-_\mathrm{max}}\int\limits_{\{u=t\}\cap N}\sqrt{\frac{W}{W_R}}\,d\s\ge|\Sigma_N|.$$
\end{itemize}
\end{lemma}

\proof \textit{Part I)} Observe that 
$$u_\mathrm{max}=\frac{1}{2}-\frac{n}{2(n-2)}R^2\qquad\text{and}\qquad u=\frac{1}{2}-\frac{\p^2}{2}-\frac{R^n}{(n-2)\p^{n-2}}$$
hence we have
$$u_\mathrm{max}-u=\frac{\p^2}{2}-\frac{n}{2(n-2)}R^2+\frac{R^n}{(n-2)\p^{n-2}},$$
and, using \eqref{W_R}, 
$$\frac{W_R}{u_\mathrm{max}-u}=\frac{2(n-2)(\p^n-R^n)^2}{\p^n(n\p^n-2\p^n-nR^2\p^{n-2}+2R^n)}.$$
Denoting by $z:=\p^n-R^n$, as $x\to p$ with $x\in N$ then $z\to0$ and we get that
\begin{equation}\label{umax-u}
u_\mathrm{max}-u=\frac{1}{2n R^{2n}}z^2+\mathcal{O}(z^3)
\end{equation}
and
\begin{equation}\label{WR/umax-u}
\frac{W_R}{u_\mathrm{max}-u}=2n-\frac{4(n-1)}{3R^n}z+\mathcal{O}(z^2), 
\end{equation}
which proves \eqref{lim WR/u-umax}. From \eqref{umax-u}, we get that
$$z=\begin{cases} 
R^{n-1}\sqrt{2n(u_\mathrm{max}-u)}+\mathcal{O}(u_\mathrm{max}-u)\hfill\text{if}\,\,\,\tau(u,N)<\sqrt{n}\\
-R^{n-1}\sqrt{2n(u_\mathrm{max}-u)}+\mathcal{O}(u_\mathrm{max}-u)\qquad\text{if}\,\,\,\tau(u,N)>\sqrt{n}\end{cases}$$
and plugging into \eqref{WR/umax-u}, we obtain \eqref{stima W_R}.\\
\textit{Part II)} We recall that from \cite[Theorem 3.1]{BoChMa} for every $p\in\Sigma$ it holds
\begin{equation}\label{umax-u BoChMa}
u=u_\mathrm{max}-\frac{n}{2}r^2-\frac{n}{6}\HH r^3+\mathcal{O}_1(r^4)
\end{equation}
on a certain $V$ neighborhood of $p$ such that $V\cap\mathrm{MAX}(u)\subset\Sigma$ and $V\setminus\Sigma$ as two connected components. 
The expansion for $W$ comes from \eqref{umax-u BoChMa}, by computing the squared norm of its gradient. The expansion for $W_R$ follows plugging expansion \eqref{umax-u BoChMa} into \eqref{stima W_R}.\\
\textit{Part III)} From the previous estimates $\frac{W}{W_R}\to0$ when $t\to u^-_\mathrm{max}$, so to prove the thesis, it is sufficient to show that
$$\limsup\limits_{t\to u^-_\mathrm{max}}\big|\{u=t\}\cap \overline{N}\big|\ge|\Sigma_N|.$$
The proof of this inequality is the same as \eqref{limsup>0}.
\endproof

%----------------------------------------------------------------------------

\section{Mean curvature bounds}

In this section, we will see how the gradient estimates of Theorem \ref{stime gradiente} can be fruitfully exploited to obtain bounds on the mean curvature of the boundary components of a connected component $N$ of $\Omen$. 

First, we obtain some a priori bounds for the mean curvature of the boundary component $\partial\Om\cap\overline{N}$ at points where the norm of the gradient restricted to $\partial\Om\cap\overline{N}$ reaches its maximum value.  

\begin{prop}\label{prop H bound 1}
Let $(\Om,u)$ be a solution to problem \eqref{1} and $N$ be a connected component of $\Omen$. Let $R=R(u,N)\in[0, R_\mathrm{max})$ be the expected core radius of $u$ in the region $N$, $r_1$ and $r_2$ be the two positive zeros of the function defined in \eqref{f_R}, and assume that Normalisation 1 holds. 

Then, for every $p\in\partial\Om\cap\overline{N}$ such that 
$$|\nabla u|(p)=\max_{\partial\Om\cap\overline{N}}|\nabla u|,$$
we have that
$$\HH(p)\le\begin{cases}
\frac{n-1}{r_2}\qquad\text{if }\tau(u,N)<\sqrt{n}\\
-\frac{n-1}{r_1}\hfill\text{if }\tau(u,N)>\sqrt{n}
\end{cases}$$
where $\HH(p)$ is the mean curvature of $\partial\Om\cap\overline{N}$ at the point $p$, computed with respect to the unit normal pointing outside $N$.
\end{prop}

\proof
We compute the Taylor expansion of both $W$ and $W_R$ along the segment $\gamma(s)=p+\frac{\nabla u}{|\nabla u|}(p) s$. For what concerns $W$, we have
$$W\left(\gamma(s)\right)=W(p)+2\nabla^2u(\nabla u,\nabla u)\frac{1}{|\nabla u|}s+o(s).$$
Recalling that $\HH=\frac{\nabla^2 u(\nabla u, \nabla u)+n|\nabla u|^2}{|\nabla u|^3}$, so that $\nabla^2u(\nabla u,\nabla u)=(H\sqrt{W}-n)W$, we get
$$W\left(\gamma(s)\right)=W(p)-2\Big(n-\HH\sqrt{W(p)}\Big)\sqrt{W(p)}s+o(s).$$
On the other hand, recalling that $\nabla W_R=-2\frac{\p^n+(n-1)R^n}{\p^n}\nabla u$, we obtain
$$W_R(\gamma(s))=W_R(p)-2\frac{\p^n+(n-1)R^n}{\p^n}\sqrt{W}s+o(s).$$
Since by assumption $W(p)=W_R(p)$, and $W\le W_R$ on $N$ from Theorem \ref{stime gradiente}, we finally have
$$\HH(p)\le (n-1)\frac{(\p^n-R^n)}{\p|\p^n-R^n|}.$$
We conclude by observing that if $\tau(u,N)<\sqrt{n}$, then $\p(p)=\psi_+(0)=r_2$, while if $\tau(u,N)>\sqrt{n}$, then $\p(p)=\psi_-(0)=r_1$.
\endproof

Now we prove the analogous mean curvature bounds for $\Sigma_N$. 

\begin{prop}\label{prop 4.2 new}
Let $(\Om,u)$ be a solution to problem \eqref{1} and $N$ be a connected component of $\Omen$. Let $R=R(u,N)\in[0, R_\mathrm{max}]$ be the expected core radius of $u$ in the region $N$, assume that Normalisation 1 holds, and suppose that $\Sigma_N\not=\emptyset$

Then, for every $p\in \Sigma_N$ it holds
$$\HH(p)\le\begin{cases}
-\frac{n-1}{R}\quad&\text{if }\tau(u,N)<\sqrt{n}\\
\frac{n-1}{R}&\text{if }\tau(u,N)\ge\sqrt{n}
\end{cases}$$
where $\HH(p)$ is the mean curvature of $\Sigma_N$ at the point $p$, computed with respect to the unit normal pointing outside $N$.
\end{prop}

\proof
We take advantage of the Taylor expansion for $W$ and $W_R$ stated in Lemma \ref{stime appendix} II) and Theorem \ref{stime gradiente}, to obtain
$$\HH\le \begin{cases}
\frac{\HH}{3}-\frac{2n-2}{3R}\quad&\text{if }\tau(u,N)<\sqrt{n}\\
\frac{\HH}{3}+\frac{2n-2}{3R}&\text{if }\tau(u,N)>\sqrt{n}.
\end{cases}$$
When $\tau(u,N)\not =\sqrt{n}$, the desired estimates follows at once. 
The case $\tau(u,N)=\sqrt{n}$, which corresponds to $R(u,N)=R_{\rm max}$, can be treated in the following way: we 
rescale $u$ in order to have 
$u_{\rm max}=(u_R)_{\rm max}$,
with $R=R_\mathrm{max}-\e$ and 
observe that the very same argument of
the proof of Theorem \ref{stime gradiente} 
yields $W\leq W_R$ in $N$.
We can then use Lemma \ref{stime appendix} II) as above and finally conclude, by passing to the limit as $\e\to0^+$.
\endproof

%----------------------------------------------------------------------------

\section{Proof of Theorem \ref{theo d new} and Theorem \ref{theo b new}}

First, we employ the previous a priori mean curvature bounds to prove a \textit{weaker} version of Theorem \ref{theo d new}.

\begin{theorem}\label{theo d partial}
Let $(\Om,u)$ be a solution to problem \eqref{1}, such that $\Om$ is a ring-shaped domain whose inner and outer boundary components are denoted by $\Gamma_1$ and $\Gamma_2$ respectively. Assume that there exists a smooth hypersurface $\Sigma\subseteq\mathrm{MAX}(u)$ dividing $\Om$ into two regions $\Om_1$ and $\Om_2$, such that $\partial\Om\cap\overline{\Om}_1 =\Gamma_1$ and $\partial\Om\cap\overline{\Om}_2= \Gamma_2$. In particular, $R_1:=R(u,\Gamma_1)$ and $R_2:=R(u,\Gamma_2)$ are well-defined. Finally, suppose that $\tau(u,\Om_2)<\sqrt{n}$. 

Then, for every $p\in\Sigma$ we have
\begin{equation}\label{stima <H<}
\frac{n-1}{R_2}\sqrt{(u_{R_2})_\mathrm{max}}\le \HH(p)\sqrt{u_\mathrm{max}} \le \frac{n-1}{R_1}\sqrt{(u_{R_1})_\mathrm{max}}
\end{equation}
where $\HH(p)$ mean curvature of $\Sigma$ with respect to the exterior unit normal to $\Om_1$. Moreover, $$R_2\ge R_1$$ 
and $R_1=R_2=:R$ if and only if $(\Om,u)$ is equivalent to $(\Om_R,u_R)$.
\end{theorem}

\proof
Since $\Omega$ is divided into two regions $\Omega_1$ and $\Om_2$ such that  $\partial\Om\cap\overline{\Om}_i=\Gamma_i$
for $i=1,2$, then $\tau(u,\Gamma_i)=\tau(u,\Omega_i)$. From Theorem \ref{teo2.1new} we know that $\tau(u,\Omega_i)\ge1$, it follows immediately that the expected core radius associated with the two boundary components $\Gamma_1$ and $\Gamma_2$ is well-defined. 

Let us consider the region $\Omega_1$. We normalize the solution as in \eqref{scaling inv} taking $\la=\sqrt{(u_{R_1})_\mathrm{max}/u_\mathrm{max}}=:\la_1$ (and $c=0$). 
In such a way the rescaled solution $(\Om_{\la_1},u_{\la_1})$ is such that $\tau(u_{\la_1},\Om_{\la_1})=\tau(u,\Om_1)$. We apply Proposition \ref{prop 4.2 new}, where we denote by $\HH_1$ the mean curvature of $\la_1\Sigma$ with respect to the unit normal pointing outside $\Om_{\la_1}$. If $\tau(u,\Om_1)<\sqrt{n}$, it holds
$$\HH_1(p)\le-\frac{n-1}{R_1}<0\qquad \forall p\in\la_1\Sigma$$
which is not possible since any closed hypersurface in $\R^n$ admits at least one point having positive mean curvature. Hence we have $\tau(u,\Om_1)\ge\sqrt{n}$ and again by Proposition \ref{prop 4.2 new} it holds
$$\HH_1(p)\le\frac{n-1}{R_1}\qquad \forall p\in\la_1\Sigma.$$
Recalling that when applying a homothety by $a\in\R^+$ the mean curvature rescales by a factor $1/a$, we get
\begin{equation}\label{Hhh}
\frac{\HH(p)}{\la_1}\le\frac{n-1}{R_1}  
\end{equation}
where $\HH$ denotes the mean curvature of $\Sigma$  with respect to the unit normal pointing outside $\Om_1$. For what concerns $\Om_2$, since we have
\begin{equation}\label{Hhhh}
-\frac{\HH(p)}{\la_2}=\HH_2(p)\le-\frac{n-1}{R_2}   
\end{equation}
where $\HH_2$ denotes the mean curvature of $\la_2\Sigma$ with respect to the unit normal pointing outside $\Om_2$ and $\la_2=\sqrt{(u_{R_2})_\mathrm{max}/u_\mathrm{max}}$. The above inequality follows from Proposition  \ref{prop 4.2 new} coupled with the hypothesis $\tau(u,\Om_2)<\sqrt{n}$.
Putting together \eqref{Hhh} and \eqref{Hhhh} we get \eqref{stima <H<}.

Now, note that as a consequence of definition \eqref{def umax} and \eqref{stima <H<} we have $f(R_2)\le f(R_1)$ where $f:R\mapsto (n-1)\sqrt{\frac{1}{2R^2}-\frac{n}{2n-4}}$. Since the function $f$ is non-increasing for every $R\in(0,R_{\rm max})$ it follows that $R_1\le R_2$.

Finally, the rigidity part of the statement follows from the fact that if $R_1=R_2=R$ then $\HH=\frac{n-1}{R}$ on $\Sigma$, and hence $\Sigma$ is the $(n-1)$-dimensional sphere of radius $R$. We can thus conclude, by combining the uniqueness of solutions given by the Cauchy-Kovalevskaya Theorem and the analyticity of $u$, as in the proof Theorem \ref{stime gradiente}.
\endproof

\proof[Proof of Theorem \ref{theo d new}]
It is sufficient to prove that there exists an hypersurface $\Sigma\subseteq\mathrm{MAX}(u)$ dividing the set $\Om$ into two regions  $\Om_1$ and $\Om_2$, such that $\partial\Om\cap\bar{\Om}_1=\Gamma_1$ and $\partial\Om\cap\bar{\Om}_2= \Gamma_2$, so that we can conclude by exploiting Theorem \ref{theo d partial}.

Since $\mathcal{H}^{n-1}(\mathrm{MAX}(u))>0$, then $\Sigma^{n-1}$ top-stratum of $\mathrm{MAX}(u)$ is nonempty (refer to point (F2) in the Introduction for further details on the {\L}ojasiewicz Structure Theorem). Moreover, since $\nabla^2 u\not=0$ in $\mathrm{MAX}(u)$, from \cite[Corollary 3.4]{BoChMa} we get that $\Sigma^{n-1}=\overline{\Sigma^{n-1}}$ and $\Sigma^{n-1}$ is a complete real analytic hypersurface with empty boundary. Hence, the set $\mathrm{MAX}(u)$ is made up of a finite number of complete hypersurfaces with empty boundaries and a finite number of lower-dimensional strata. 
A closer look shows that the top-stratum of $\mathrm{MAX}(u)$ is given by only one closed hypersurface dividing the set $\Om$ as required. This can be proven by arguing by contradiction and using the strong maximum principle to rule out both the existence of closed hypersurfaces not enclosing the inner boundary $\Gamma_1$, and the existence of two hypersurfaces one enclosed in the other. 
\endproof

In order to prove our dichotomy theorem we make use of the following result, which is a sub-case of \cite[Theorem 2]{Sirakov}.

\begin{theorem}\label{sirakov}
Let $\Om$ be a ring-shaped domain where $\Gamma_\mathrm{in}$ and $\Gamma_\mathrm{out}$ denote the inner and the outer boundary components, respectively. Assume that $u$ solves
\begin{equation}\label{sist sirak}
\begin{cases}
    \Delta u=c\qquad&\text{in }\Om,\\
    u=a &\text{on }\Gamma_\mathrm{in},\\
    u=0   &\text{on }\Gamma_\mathrm{out},
\end{cases}
\end{equation}
where $a$ is a positive constant and $c$ a nonzero constant. Suppose that $|\nabla u|$ is locally constant on $\partial \Om$ and, denoting by $\nu$ the exterior unit normal to $\Om$, suppose that $\partial_\nu u\ge0$ on $\Gamma_\mathrm{in}$.
Then $(\Om,u)$ is rotationally symmetric.    
\end{theorem}

This theorem has been proven in \cite{Sirakov} using the method of moving planes. Another possible proof relies on the comparison geometry technique described in the previous sections. Indeed, such arguments have been used in \cite{Borghini} to prove the two dimensional version or Theorem \ref{sirakov}: we believe that they can be extended to arbitrary dimension without major modifications.

\proof[Proof of Theorem \ref{theo b new}]
Let $u$ be a solution to \eqref{1} and denote by $u_{\rm max}=\max_{\overline\Om} u>0$. As in the proof of Theorem \ref{theo d new}, one can prove that the top-stratum of $\mathrm{MAX}(u)$ is given by only one closed hypersurface $\Sigma$ dividing the set $\Om$ into two regions  $\Om_1$ and $\Om_2$, such that $\partial\Om\cap\overline{\Om}_1=\Gamma_1$ and $\partial\Om\cap\overline{\Om}_2= \Gamma_2$. Let us first consider the case when $|\nabla u|$ is constant on the outer boundary component $\Gamma_2$. In this case, $u$ solves problem \eqref{sist sirak} in
$\Om_2$, with $c=-n$, $\Gamma_\mathrm{in}=\Sigma$, 
$\Gamma_\mathrm{out}=\Gamma_2$,
and with $a=u_{\rm max}$. 
Since $\partial_\nu u=0$ on $\Sigma$, a direct application of Theorem \ref{sirakov} yields the radial symmetry of the solution in $\Om_2$. In turn, this implies the rotational symmetry on the whole $\Omega$, arguing as in end of the proof of Theorem \ref{stime gradiente}. 
Let us now consider the case when $|\nabla u|$ is constant on $\Gamma_1$. In this case, consider the function $v=u_\mathrm{max}-u$, which solves problem \eqref{sist sirak} in $\Om_1$ with $c=n$, $a=u_\mathrm{max}$, $\Gamma_\mathrm{in}=\Gamma_1$, 
$\Gamma_\mathrm{out}=\Sigma$, and with $a=u_{\rm max}$.
Moreover, $\partial_\nu v$ equals a positive constant on $\Gamma_\mathrm{in}$, due to the Hopf Lemma. Therefore, also in this case we can apply Theorem \ref{sirakov} to deduce that $v$, and in turn $u$, is rotationally symmetric in $\Om_1$. As in the previous case, one can then deduce rotational symmetry throughout $\Om$.
\endproof

%---------------------------------------------------------------------------

\section{Proof of the area bounds}

In this section, we obtain some \textit{a priori} bounds on the area of the boundary components of a given region $N\subseteq \Omen$. We distinguish between the bounds holding for the boundary portion of the region $N$ on $\partial\Om$, that is $\Gamma_N:=\partial\Om \cap\overline N$, and the ones holding for a smooth portion of the top stratum of $\bar N\cap \mathrm{MAX}(u)$, that we denote by $\Sigma_N$ (see Proposition \ref{prop bnd ar 1} and Proposition \ref{prop bnd ar 2} below). We conclude this section with the proof of the third area bound, where the area of $\Sigma_N$ is compared to the area of $\Gamma_N$ (see Theorem \ref{theo 2.4}).

\begin{prop}\label{prop bnd ar 1}
Let $(\Om,u)$ be a solution to problem \eqref{1} and $N$ be a connected component of $\Omen$. Let $R=R(u,N)\in[0, R_\mathrm{max})$ be the expected core radius of $u$ in the region $N$, and $r_1$ and $r_2$ the two positive zeros of the function \eqref{f_R}. Assume that Normalisation 1 holds. Finally, suppose that $\Gamma_N$ is connected and the $|\nabla u|$ is constant on $\Gamma_N$. 
\begin{itemize}
    \item[i)] If $\tau(u,N)<\sqrt{n}$, then 
    $$|\Gamma_N|\ge r_2\,\int_{\Gamma_N}\frac{\HH}{n-1}d\s,$$
    where $\HH$ is the mean curvature of $\Gamma_N$ computed with respect to the unit normal pointing outside $N$.
    \item[ii)] If $\tau(u,N)>\sqrt{n}$, then 
    $$|\Gamma_N|\le\left(\int_{\Gamma_N}\frac{\HH}{n-1}d\s\right)r_1, $$
    where $\HH$ is the mean curvature of $\Gamma_N$ computed with respect to the unit normal pointing inside $N$.
\end{itemize}
Moreover, if equality holds, then $(\Om,u)$ is equivalent to the ring-shaped model solution $(\Om_R,u_R)$. 
\end{prop}

\proof
We exploit the a priori bounds on the mean curvature of the boundary components of a region $N\subseteq \Omen$ obtained in Proposition \ref{prop H bound 1}. More precisely, assuming that $|\nabla u|$ is constant on $\Gamma_N$, guarantees that such curvature bounds hold all over $\Gamma_N$. The area bounds then follow at once integrating these curvature estimates on $\Gamma_N$.
For what concerns the rigidity, if we are in the case $i)$ and equality holds in the corresponding area bound, then again using Proposition \ref{prop H bound 1} we get that $\HH=\frac{n-1}{r_2}$ all over $\Gamma_N$, one then concludes rotational symmetry via Alexandrov's Theorem. The case $ii)$ is analogous.
\endproof

\begin{obs} 
Under the hypothesis of the previous theorem, we note that in the case $\tau(u,N)<\sqrt{n}$, we can say something more, using the Minkowsky inequality. Indeed, if $\Gamma_N$ encloses an outward minimizing 
set, then the Minkowski inequality reads 
$$\int_{\Gamma_N}\frac{H}{n-1}\,d\s\ge |\mathbb{S}^{n-1}|^\frac{1}{n-1} |\Gamma_N|^\frac{n-2}{n-1},$$
where $\HH$ is the mean curvature of $\Gamma_N$ computed with respect to the unit normal pointing outside $N$. 
This inequality, coupled with point $i)$ in the previous proposition, 
yields 
$$r_2\le\left(\frac{|\Gamma_N|}{|\mathbb{S}^{n-1}|}\right)^\frac{1}{n-1}.$$
Moreover, the equality holds if and only if the solution is rotationally symmetric.
\end{obs}

\begin{prop}\label{prop bnd ar 2}
Let $(\Om,u)$ be a solution to problem \eqref{1} and $N$ be a connected component of $\Omen$. Let $R=R(u,N)\in[0, R_\mathrm{max})$ be the expected core radius of $u$ in the region $N$, and assume that Normalisation 1 holds. 
\begin{itemize}
    \item[(i)] If $\tau(u,N)<\sqrt{n}$, then $$|\Sigma_N|\le\left(\int_{\Sigma_N}\frac{\HH}{n-1}\right)R$$
    where $\HH$ is the mean curvature of $\Sigma_N$ with respect to the unit normal pointing inside $N$.
    \item[(ii)] If $\tau(u,N)\ge\sqrt{n}$, then 
    $$|\Sigma_N|\ge\left(\int_{\Sigma_N}\frac{\HH}{n-1}\right)R$$
    where $\HH$ is the mean curvature of $\Sigma_N$ with respect to the unit normal pointing outside $N$.
\end{itemize}
Moreover, if equality holds, then $(\Om,u)$ is equivalent to the ring-shaped model solution $(\Om_R,u_R)$. 
\end{prop}

\proof 
Since $\Sigma_N$ is the top-stratum of $\mathrm{MAX}(u)\cap\overline{N}$, then we can take advantage of the curvature bounds obtained in Proposition \ref{prop 4.2 new}. Namely, for every $p\in \Sigma_N$ it holds
$$\HH(p)\le\begin{cases}
\frac{n-1}{R}\qquad\text{if }\tau(u,N)<\sqrt{n}\\
-\frac{n-1}{R}\hfill\text{if }\tau(u,N)\ge\sqrt{n}
\end{cases}$$
where $\HH(p)$ is the mean curvature of $\Sigma$ at $p$, computed with respect to the unit normal pointing outside $N$. Integrating these curvature bounds on $\Sigma_N$ allows us to conclude. As in the proof of Proposition \ref{prop bnd ar 1} the rigidity part can be obtained as a consequence of Alexandrov Theorem.
\endproof

\proof[Proof of Proposition \ref{theo 2.4}]
First of all, to prevent blow-up on the considered quantities on $\Sigma_N$, we set $N_\e=N\cap\{u\le u_\mathrm{max}-\e\}$. By Sard's Theorem the regular values of $u$ are dense, so we can find a sequence $(\e_i)_{i\in\mathbb{N}}$ of positive real numbers converging to zero such that for every $i$ the level set $\{u=u_\mathrm{max}-\e_i\}$ is regular. Using the pseudo-radial function introduced in Definition \ref{def Phi}, note first that 
$$\int\limits_{N_\e}\frac{\Delta u}{\p^n-R^n}dx=-\int\limits_{N_\e} \langle\nabla\left(\frac{1}{\p^n-R^n}\right),\nabla u\rangle \, dx+\int\limits_{\partial N_\e}\frac{1}{\p^n-R^n}\,\langle\nabla u,\nu\rangle\,d\s,$$
where $\nu$ is the outward unit normal to $\partial N_\e$. Recalling that $u$ solves \eqref{1} and splitting $\partial N_\e$ into $\Gamma_N$ and $\Sigma_{N,\e}:=N\cap\{u=u_\mathrm{max}-\e\}$, we get
$$\int\limits_{N_\e}\frac{-n}{\p^n-R^n}+\frac{n\p^{2n-2}}{(\p^n-R^n)^3}|\nabla u|^2dx=\int\limits_{\Sigma_{N,\e}}\frac{|\nabla u|}{\p^n-R^n}\,d\s-\int\limits_{\Gamma_N}\frac{|\nabla u|}{\p^n-R^n}\,d\s,$$
because $\nu=-\nabla u/|\nabla u|$ on $\Gamma_N$ and $\nu=\nabla u/|\nabla u|$ on $\Sigma_{N,\e}$. By definition of $W$ and by \eqref{W_R}, we have thus obtained that
\begin{equation}\label{integral id area bnds}
\int\limits_{N_\e}\frac{n\p^{2n-2}}{(\p^n-R^n)^3}(W-W_R)\,dx=\int\limits_ {\Sigma_{N,\e}} \frac{|\nabla u|}{\p^n-R^n}\,d\s-\int\limits_{\Gamma_N} \frac{|\nabla u|}{\p^n-R^n}\,d\s.    
\end{equation}
We finally derive the desired bounds from this identity, exploiting the fact that $W-W_R\leq0$, due to Theorem \ref{stime gradiente}.
At the same time, the behavior of the pseudo-radial function $\p$ depends on whether $\tau(u,N)<\sqrt{n}$ or $\tau(u,N)>\sqrt{n}$.

\textbf{Case $\tau(u,N)<\sqrt{n}$.} By definition $\p=\p_+\in(R,r_2]$ over $\overline N_\e$, so that the left-hand-side of \eqref{integral id area bnds} is non-positive.
Moreover, $\p_+\equiv r_2$ on $\Gamma_N$ and $\max_{\Gamma_N}|\nabla u|=\frac{r_2^n-R^n}{r_2^{n-1}}$
(see Remark \ref{spegazione W WR}), so that 
$$\int\limits_{\Gamma_N} \frac{|\nabla u|}{\p_+^n-R^n}\,d\s = \int\limits_{\Gamma_N} \frac{|\nabla u|}{r_2^n-R^n}\,d\s= \int\limits_{\Gamma_N} \frac{1}{r_2^{n-1}}\,\,\frac{|\nabla u|}{\max_{\Gamma_N}|\nabla u|}\,d\s\le\,\frac{|\Gamma_N|}{r_2^{n-1}}.$$
We compute the limit of the integral over $\Sigma_{N,\e}$ as $\e\to0^+$
$$\lim_{\e\to0^+}\int\limits_ {\Sigma_{N,\e}} \frac{|\nabla u|}{\p_+^n-R^n}\,d\s= \lim_{\e\to0^+}\int\limits_ {\Sigma_{N,\e}}\frac{1}{\p_+^{n-1}} \sqrt{\frac{W}{W_R}}\,d\s= \frac{1}{R^{n-1}}\lim_{\e\to0^+}\int\limits_ {\Sigma_{N,\e}}\sqrt{\frac{W}{W_R}}\,d\s\ge\frac{|\Sigma_N|}{R^{n-1}},$$
where in the last inequality we used point $III)$ in Lemma \ref{stime appendix}. Combining all these facts, we get
\begin{equation}
0\ge \lim_{\e\to0^+} \int\limits_ {\Sigma_{N,\e}} \frac{|\nabla u|}{\p_+^n-R^n}\,d\s -\int\limits_{\Gamma_N} \frac{|\nabla u|}{\p_+^n-R^n}\,d\s\ge \frac{|\Sigma_N|}{R^{n-1}}-\frac{|\Gamma_N|}{r_2^{n-1}},
\end{equation}
which proves point $i)$.
Now, if $\frac{|\Sigma_N|}{R^{n-1}}=\frac{|\Gamma_N|}{r_2^{n-1}}$, then
from \eqref{integral id area bnds} by monotone convergence we get that 
$$\int\limits_{N}\frac{n\p^{2n-2}}{(\p^n-R^n)^3}(W-W_R)\,dx=0$$
this fact couple with Theorem \eqref{stime gradiente} gives that $W=W_R$ in the whole $N$. We then conclude that $(\Om,u)$ is a ring-shaped model solution with core radius $R$, by exploiting the rigidity part of the same theorem.

\textbf{Case $\tau(u,N)>\sqrt{n}$.} Follows the same argument as in the previous case. We only recall that when $\tau(u,N)>\sqrt{n}$ then $\p=\p_-\in[r_1,R)$ over $\overline N_\e$, so the left-hand-side of \eqref{integral id area bnds} is non-negative. Moreover, $\p_-\equiv r_1$ on $\Gamma_N$ and $\max_{\Gamma_N}|\nabla u|=\frac{R^n-r_1^n}{r_1^{n-1}}$, so 
$$\lim_{\e\to0^+}\int\limits_ {\Sigma_{N,\e}} \frac{|\nabla u|}{\p_-^n-R^n}\,d\s\le-\frac{|\Sigma_N|}{R^{n-1}}.$$
\endproof

\begin{ackn}The authors are
  members of the Gruppo Nazionale per l'Analisi Matematica, la
  Probabilit\`a e le loro Applicazioni (GNAMPA), which is part of the
  Istituto Nazionale di Alta Matematica (INdAM).
  They gratefully acknowledge the support 
  by the PRIN project ``Contemporary perspectives on geometry and gravity''.
\end{ackn}

%-----------------------------------------------------------------------------

\bibliography{main}

\begin{thebibliography}{10}

\bibitem{AgoBorMaz}
{\sc Agostiniani, V., Borghini, S., and Mazzieri, L.}
\newblock On the serrin problem for ring-shaped domains.
\newblock {\em To appear in Journal of the European Mathematical Society\/} (2023).

\bibitem{Alex}
{\sc Alexandrov, A.~D.}
\newblock A characteristic property of spheres.
\newblock {\em Ann. Mat. Pura Appl. (4) 58\/} (1962), 303--315.

\bibitem{Borghini}
{\sc Borghini, S.}
\newblock Symmetry results for {S}errin-type problems in doubly connected domains.
\newblock {\em Math. Eng. 5}, 2 (2023), Paper No. 027, 16.

\bibitem{BoChMa}
{\sc Borghini, S., Chru\'{s}ciel, P.~T., and Mazzieri, L.}
\newblock On the uniqueness of {S}chwarzschild--de {S}itter spacetime.
\newblock {\em Arch. Ration. Mech. Anal. 247}, 2 (2023), Paper No. 22, 35.

\bibitem{BorMazCQG}
{\sc Borghini, S., and Mazzieri, L.}
\newblock On the mass of static metrics with positive cosmological constant: {I}.
\newblock {\em Classical Quantum Gravity 35}, 12 (2018), 125001, 43.

\bibitem{Dam}
{\sc Damascelli, L., Pacella, F., and Ramaswamy, M.}
\newblock Symmetry of ground states of {{\(p\)}}-{Laplace} equations via the moving plane method.
\newblock {\em Arch. Ration. Mech. Anal. 148}, 4 (1999), 291--308.

\bibitem{Es_Ma}
{\sc Espinar, J.~M., and Marin, D.~A.}
\newblock An overdetermined eigenvalue problem and the critical catenoid conjecture.
\newblock ArXiv Preprint Server {--} arXiv:2310.06705.

\bibitem{Fall}
{\sc Fall, M.~M., and Jarohs, S.}
\newblock Overdetermined problems with fractional {Laplacian}.
\newblock {\em ESAIM, Control Optim. Calc. Var. 21}, 4 (2015), 924--938.

\bibitem{GiNiNi}
{\sc Gidas, B., Ni, W.-M., and Nirenberg, L.}
\newblock Symmetry and related properties via the maximum principle.
\newblock {\em Commun. Math. Phys. 68\/} (1979), 209--243.

\bibitem{KamSci}
{\sc Kamburov, N., and Sciaraffia, L.}
\newblock Nontrivial solutions to {Serrin}'s problem in annular domains.
\newblock {\em Ann. Inst. Henri Poincar{\'e}, Anal. Non Lin{\'e}aire 38}, 1 (2021), 1--22.

\bibitem{KraPar}
{\sc Krantz, S.~G., and Parks, H.~R.}
\newblock {\em A primer of real analytic functions}, second~ed.
\newblock Birkh\"{a}user Advanced Texts: Basler Lehrb\"{u}cher. [Birkh\"{a}user Advanced Texts: Basel Textbooks]. Birkh\"{a}user Boston, Inc., Boston, MA, 2002.

\bibitem{Morrey}
{\sc Morrey, Jr., C.~B.}
\newblock On the analyticity of the solutions of analytic non-linear elliptic systems of partial differential equations. {I}. analyticity in the interior.
\newblock {\em Amer. J. Math. 80\/} (1958), 198--218.

\bibitem{Reichel}
{\sc Reichel, W.}
\newblock Radial symmetry by moving planes for semilinear elliptic {BVPs} on annuli and other non-convex domains.
\newblock In {\em Elliptic and parabolic problems. Proceedings of the 2nd European conference, Pont-\`a-Mousson, June 1994}. Harlow: Longman Scientific \& Technical; New York, NY: Wiley, 1995, pp.~164--182.

\bibitem{ReichelARMA}
{\sc Reichel, W.}
\newblock Radial symmetry for elliptic boundary-value problems on exterior domains.
\newblock {\em Arch. Ration. Mech. Anal. 137}, 4 (1997), 381--394.

\bibitem{Serrin}
{\sc Serrin, J.}
\newblock A symmetry problem in potential theory.
\newblock {\em Arch. Rational Mech. Anal. 43\/} (1971), 304--318.

\bibitem{SerrinZou}
{\sc Serrin, J., and Zou, H.}
\newblock Symmetry of ground states of quasilinear elliptic equations.
\newblock {\em Arch. Ration. Mech. Anal. 148}, 4 (1999), 265--290.

\bibitem{Sirakov}
{\sc Sirakov, B.}
\newblock Symmetry for exterior elliptic problems and two conjectures in potential theory.
\newblock {\em Annales de l'Institut Henri Poincare (C) Analyse Non Lineaire 18}, 2 (2001), 135 – 156.

\bibitem{SouSou}
{\sc Sou\v{c}ek, J., and Sou\v{c}ek, V.}
\newblock Morse-{S}ard theorem for real-analytic functions.
\newblock {\em Comment. Math. Univ. Carolinae 13\/} (1972), 45--51.

\end{thebibliography}
\bibliographystyle{acm} %abbrv}	

\end{document}